\documentclass[11pt]{article}
\usepackage{amsmath,amsthm,amssymb}
\usepackage{fancyhdr}
\usepackage[british]{babel}
\usepackage{geometry,mathtools}
\usepackage{enumitem}
\usepackage{algpseudocode}
\usepackage{dsfont}
\usepackage{centernot}
\usepackage{xstring}
\usepackage{colortbl}
\usepackage[section]{algorithm}
\usepackage{graphicx}
\usepackage[font={footnotesize}]{caption}
\usepackage[usenames,dvipsnames,table]{xcolor}
\usepackage{pstricks,tikz}
\usepackage[h]{esvect}
\usepackage[
		bookmarksopen=true,
		bookmarksopenlevel=1,
		colorlinks=true,
		linkcolor=darkblue,
        linktoc=page,
		citecolor=darkblue,
]{hyperref}
\usepackage{bbm}

\numberwithin{equation}{section}

\definecolor{darkblue}{rgb}{0,0,0.5}

\newdimen\margin
\def\textno#1&#2\par{
   \margin=\hsize
   \advance\margin by -4\parindent
          \setbox1=\hbox{\sl#1}
   \ifdim\wd1 < \margin
      $$\box1\eqno#2$$
   \else
      \bigbreak
      \hbox to \hsize{\indent$\vcenter{\advance\hsize by -3\parindent
      \it\noindent#1}\hfil#2$}
      \bigbreak
   \fi}

\newtheorem{theorem}[algorithm]{Theorem}
\newtheorem{prop}[algorithm]{Proposition}
\newtheorem{lemma}[algorithm]{Lemma}

\theoremstyle{definition}

\newtheorem{conj}[algorithm]{Conjecture}

\newtheorem{remark}[algorithm]{Remark}

\def\proof{\removelastskip\penalty55\medskip\noindent\begin{stepenv}\end{stepenv}{\bf Proof. }} 

\def\lateproof#1{\removelastskip\penalty55\medskip\noindent\begin{stepenv}\end{stepenv}{\bf Proof of #1. }} 

\def\noproof{{\unskip\nobreak\hfill\penalty50\hskip2em\hbox{}\nobreak\hfill%
       $\square$\parfillskip=0pt\finalhyphendemerits=0\par}\goodbreak}
\def\endproof{\noproof\bigskip}

\def\claimproof{\removelastskip\penalty55\medskip\noindent{\em Proof of claim: }}
\def\noclaimproof{{\unskip\nobreak\hfill\penalty50\hskip2em\hbox{}\nobreak\hfill%
       $\square$\parfillskip=0pt\finalhyphendemerits=0\par}\goodbreak}

\def\endclaimproof{\noclaimproof\medskip}

\newcounter{stepenv}
\newenvironment{stepenv}[1][]{\refstepcounter{stepenv}}{}

\newcounter{step}[stepenv]

\newcounter{substep}[step]
\renewcommand{\thesubstep}{\thestep.\arabic{substep}}

\newcounter{claim}[stepenv]
\newenvironment{claim}[1][]{\refstepcounter{claim}\par\medskip\noindent%
        \textit{Claim~\theclaim. #1} \itshape\rmfamily}{\medskip}

\newcommand{\cE}{\mathcal{E}}

\newcommand{\cO}{\mathcal{O}}

\newcommand{\bN}{\mathbb{N}}

\newcommand{\bR}{\mathbb{R}}

\def\eps{{\epsilon}}
\newcommand{\eul}{{e}}

\newcommand{\defn}{\emph}

\newcommand{\prob}[1]{\mathrm{\mathbb{P}}\left[#1\right]}
\newcommand{\cprob}[2]{\prob{#1 \;\middle|\; #2}}
\newcommand{\expn}[1]{\mathrm{\mathbb{E}}\left[#1\right]}
\newcommand{\cexpn}[2]{\expn{#1 \;\middle|\; #2}}

\def\sm{\setminus}
\newcommand{\Set}[1]{\{#1\}}
\newcommand{\set}[2]{\{#1\,:\;#2\}}
\def\In{\subset}

\newcommand{\IND}{\mathbbm{1}}

\DeclareMathOperator{\ex}{ex}

\def\COMMENT#1{}
\def\TASK#1{}
\let\TASK=\footnote             

\begin{document}

\title{An average degree condition for independent transversals}

\author{Stefan Glock \thanks{Institute for Theoretical Studies, ETH, 8092 Z\"urich, Switzerland.
Email: \href{mailto:st3fan.g10ck@gmail.com}{\nolinkurl{st3fan.g10ck@gmail.com}}.
Research supported by Dr. Max R\"ossler, the Walter Haefner Foundation and the ETH Z\"urich Foundation.}
\and Benny Sudakov \thanks{Department of Mathematics, ETH, 8092 Z\"urich, Switzerland. Email:
\href{mailto:benny.sudakov@gmail.com} {\nolinkurl{benny.sudakov@gmail.com}}.
Research supported in part by SNSF grant 200021\textunderscore196965.}}

\date{}

\maketitle

\begin{abstract} 
In 1994, Erd\H{o}s, Gy\'{a}rf\'{a}s and {\L}uczak posed the following problem: given disjoint vertex sets $V_1,\dots,V_n$ of size~$k$, with exactly one edge between any pair $V_i,V_j$, how large can $n$ be such that there will always be an independent transversal? They showed that the maximal $n$ is at most $(1+o(1))k^2$, by providing an explicit construction with these parameters and no independent transversal. They also proved a lower bound which is smaller by a $2\eul$-factor. 

In this paper, we solve this problem by showing that their upper bound construction is best possible: if $n\le (1-o(1))k^2$, there will always be an independent transversal.
In fact, this result is a very special case of a much more general theorem which concerns independent transversals in arbitrary partite graphs that are `locally sparse', meaning that the maximum degree between each pair of parts is relatively small. In this setting, Loh and Sudakov provided a global \emph{maximum} degree condition for the existence of an independent transversal. We show that this can be relaxed to an \emph{average} degree condition.

We can also use our new theorem to establish tight bounds for a more general version of the Erd\H{o}s--Gy\'{a}rf\'{a}s--{\L}uczak problem and solve a conjecture of Yuster from 1997. This exploits a connection to the Tur\'an numbers of complete bipartite graphs, which might be of independent interest.
\end{abstract}

\section{Introduction}

Given a graph $G$ and a partition $V(G)=V_1\cup \dots \cup V_m$ of its vertex set, an \defn{independent transversal} of $G$ (with respect to $\Set{V_i}_{i\in [m]}$) is an independent set in $G$ which contains exactly one vertex from each part~$V_i$. The problem of finding sufficient conditions for the existence of independent transversals was raised in the early 70s, see e.g.~\cite{BES:75,erdos:72}. It has since been extensively studied~\cite{ABZ:07,alon:88,alon:92,haxell:01,haxell:04,jin:92,LohSudakov:07,meshulam:03,reed:99,RS:02,ST:06,yuster:97,yuster:97b}, not least because the basic concept appears in many different contexts, such as SAT, linear arboricity, strong chromatic number and list colouring.

Most prominently, sufficient conditions for the existence of an independent transversal have been formulated in terms of the maximum degree $\Delta$ of $G$ and the minimum size of the parts in the partition. Bollob\'as, Erd\H{o}s and Szemer\'edi~\cite{BES:75} conjectured that if all the parts have size at least $2 \Delta$, then an independent transversal exists. An application of the Lov\'asz Local Lemma, due to Alon~\cite{alon:88}, shows that if all the parts have size at least $2\eul \Delta$, then choosing one vertex from each part at random produces an independent set with positive probability, and thus the existence of an independent transversal is guaranteed.
In a breakthrough, Haxell~\cite{haxell:01} improved the constant from $2\eul$ to $2$ using a topological proof technique, thus confirming the conjecture of Bollob\'as, Erd\H{o}s and Szemer\'edi~\cite{BES:75}. To appreciate how strong this result is, we remark that for each $\Delta$, there are constructions of graphs with maximum degree at most $\Delta$ and parts of size $2\Delta-1$ which have no independent transversal~\cite{ST:06}.
We refer to the survey~\cite{haxell:16} for more background on these results and applications.

\subsection{Locally sparse graphs}
Haxell's result yielded improvements for many combinatorial problems that can be rephrased with independent transversals. One prominent such example is list colouring, which in fact motivated much of the work on independent transversals. Given a graph $G$ with a list of colours $L_v$ for each vertex $v\in V(G)$, we aim to properly colour $G$ such that each vertex receives a colour from its list. Define the \defn{vertex-colour degree} $d_c(v)$ for a vertex $v$ and a colour $c\in L_v$ as the number of neighbours of $v$ which also have $c$ in their list. Reed~\cite{reed:99} conjectured that if $d_c(v)\le d$ for all $v\in V(G)$, $c\in L_v$, and $|L_v|\ge d+1$ for all $v\in V(G)$, then a list colouring exists. Clearly, this is true for $d=\Delta(G)$ since we can colour greedily, and the question is whether smaller lists suffice as long as the vertex-colour degrees are also small. 
This can be phrased as an independent transversal problem. For each vertex $v$, define $V_v:=\Set{v}\times L_v$, and let $\Gamma$ be the \defn{vertex-colour graph} whose vertex set is the union of the sets $V_v$, and two vertices $(v,c)$, $(v',c')$ are connected if $vv'\in E(G)$ and $c=c'$. Observe that the independent transversals of $\Gamma$ are precisely the list colourings of~$G$. Moreover, $\Delta(\Gamma)$ is the maximum vertex-colour degree. Thus, by Haxell's theorem, lists of size~$2d$ suffice. Even though Bohman and Holzman~\cite{BH:02} disproved Reed's conjecture, Reed and Sudakov~\cite{RS:02} managed to prove the conjecture asymptotically, by showing that lists of size $(1+o(1))d$ suffice.

Of course, the reason why an improvement in the constant from $2$ to asymptotically $1$ was possible in this case must lie in the special structure of the vertex-colour graph~$\Gamma$. For instance, observe that any two parts in $\Gamma$ induce a (possibly empty) matching. Aharoni and Holzman~(see~\cite{LohSudakov:07}) conjectured that this property always allows to improve the constant from $2$ to asymptotically~$1$.
This was proven in a strong form by Loh and Sudakov~\cite{LohSudakov:07}, building on earlier work of Alon~\cite{alon:92} and Reed and Sudakov~\cite{RS:02}.
Define the \defn{local degree} of $G$ (with respect to $\Set{V_i}_{i\in [m]}$) as the maximum of $\Delta(G[V_i,V_j])$ over all distinct $i,j\in [m]$.

\begin{theorem}[Loh--Sudakov~\cite{LohSudakov:07}]\label{thm:LS}
For every $\eps>0$ there exists $\gamma>0$ such that the following holds. If $G$ is a graph with $\Delta(G)\le \Delta$ whose vertex set is partitioned into parts of size at least $(1+\eps)\Delta$, and the local degree is at most $\gamma \Delta$, then $G$ has an independent transversal.
\end{theorem}

Note that the constant $1$ in the statement is optimal. For instance, let $G$ be the vertex-disjoint union of $\Delta$ cliques of order $\Delta+1$, and let each part $V_i$ consist of precisely one vertex from each clique. The local degree is $1$, the maximum degree is $\Delta$ and each part has size $\Delta$, but there is no independent transversal.

One of the main contributions of this paper is that we can replace the \emph{maximum} degree condition with an \emph{average} degree condition.
This makes the result far more applicable. To state our main theorem, define for a graph $G$ and $U\In V(G)$ the \defn{average degree} of $U$ in $G$ as $\bar{d}_G(U):=\frac{1}{|U|}\sum_{u\in U}d_G(u)$.

\begin{theorem}\label{thm:main general}
For all $\eps>0$ there exists $\gamma>0$ such that the following holds. 
Let $G$ be a graph with vertex partition $V(G)=V_1 \cup\dots\cup V_m$ and local degree at most $\gamma D$. Suppose that $|V_i|\ge (1+\eps)D$ and $\bar{d}_G(V_i)\le D$ for each $i\in[m]$. Then there exists an independent transversal.
\end{theorem}

Our main application of this new theorem is to solve an extremal problem about independent transversals by Erd\H{o}s, Gy\'{a}rf\'{a}s and {\L}uczak (see Section~\ref{sec:extremal}). A further application that we find noteworthy concerns the list colouring problem discussed earlier. The following result generalizes the asymptotic version of Reed's conjecture proved by Reed and Sudakov. Instead of requiring all vertex-colour degrees to be small, it suffices if for each vertex, the average of its colour degrees is small.

\begin{theorem}
For all $\eps>0$, the following holds for sufficiently large~$D$.
Let $G$ be a graph with given lists $\Set{L_v}_{v\in V(G)}$. If for each vertex $v$, we have $\frac{1}{|L_v|}\sum_{c\in L_v}d_c(v)\le D$ and $|L_v|\ge (1+\eps)D$, then there is a proper list colouring.
\end{theorem}

\proof
Consider the vertex-colour graph $\Gamma$ defined in the beginning of this section. The local degree is $1$, each part $V_v$ has size $|L_v|\ge (1+\eps)D$ and average degree $\frac{1}{|L_v|}\sum_{c\in L_v}d_c(v)\le D$. By Theorem~\ref{thm:main general}, $\Gamma$ has an independent transversal, which corresponds to a proper list colouring of~$G$.
\endproof

Finally, we note that asymptotically tight results on $K_t$-free transversals can be easily deduced as well. 
A \defn{$K_t$-free transversal} in a multipartite graph is a set containing exactly one vertex from each part which does not induce a~$K_t$.
Thus, a $K_2$-free transversal is precisely an independent transversal. In~\cite{LohSudakov:07}, a general result on $K_t$-free transversals, for fixed $t$, is reduced to Theorem~\ref{thm:LS}. Namely, with all other conditions unchanged, part sizes of at least $(1+o(1))\frac{\Delta}{t}$ suffice to guarantee a $K_{t+1}$-free transversal, and this is asymptotically tight.
The same reduction also works in our average degree setting.

\begin{theorem}\label{thm:clique free}
For all $\eps>0$ and $t\in \bN$ there exists $\gamma>0$ such that the following holds. 
Let $G$ be a graph with vertex partition $V(G)=V_1 \cup\dots\cup V_m$ and local degree at most $ \gamma D$. Suppose that $|V_i|\ge (1+\eps)\frac{D}{t}$ and $\bar{d}_G(V_i)\le D$ for each $i\in[m]$. Then there exists a $K_{t+1}$-free transversal.
\end{theorem}

\proof
Colour $V(G)$ with $t$ colours such that the number of monochromatic edges is minimal. Let $G'$ be the spanning subgraph containing all monochromatic edges. By the minimality of the colouring, we have $d_{G'}(v)\le d_G(v)/t$ for all $v\in V(G)$. Hence, $\bar{d}_{G'}(V_i)\le D/t$ for each $i\in[m]$, so we can apply Theorem~\ref{thm:main general} to obtain an independent transversal $T$ of~$G'$. This means that $G[T]$ is properly $t$-coloured, in particular it must be $K_{t+1}$-free.
\endproof

\subsection{Extremal problems for independent transversals}\label{sec:extremal}

In 1994, Erd\H{o}s, Gy\'{a}rf\'{a}s and {\L}uczak~\cite{EGL:94} initiated the study of independent transversals in `sparse partite hypergraphs'.
Define an \defn{$(n,k)$-graph} to be an $n$-partite graph with parts of size $k$ such that between any two parts there is exactly one edge. Clearly, any $(n,k)$-graph has $\binom{n}{2}$ edges and so is sparse when $k$ tends to infinity.
Motivated by a connection with structural graph problems, Erd\H{o}s, Gy\'{a}rf\'{a}s and {\L}uczak~\cite{EGL:94} defined the function $f(k)$ as the maximum integer $n$ such that any $(n,k)$-graph has an independent transversal. They showed that 
\begin{align}
\frac{1}{2\eul}k^2 \le f(k)\le (1+o(1))k^2 \label{EGL bound graphs}
\end{align} 
and posed the problem of estimating $f(k)$ more accurately. 
The upper bound in~\eqref{EGL bound graphs} stems from an explicit construction involving affine planes, and the lower bound is a standard application of the Lov\'asz Local Lemma (see Lemma~\ref{lem:LLL app}). 
Yuster~\cite{yuster:97} improved the $(1/2\eul)$-factor in the lower bound to~$0.65$, but no other improvement has been made until now. 
Here, we solve this problem by showing that the upper bound in~\eqref{EGL bound graphs} is tight. This is an easy consequence of Theorem \ref{thm:main general}. The crucial point is that, although an $(n,k)$-graph may have \emph{maximum} degree $n-1$, the \emph{average} degree of each part is~$(n-1)/k$.

\begin{theorem} \label{thm:EGL graph}
$f(k)=(1+o(1))k^2$.
\end{theorem}

\proof
Let $G$ be any $(n,k)$-graph. Then each part has size $k$ and average degree $(n-1)/k$, and the local degree is~$1$. Hence, if $k\ge (1+o(1))n/k$, then an independent transversal exists by Theorem~\ref{thm:main general}. This implies $f(k)\ge (1-o(1))k^2$. Together with~\eqref{EGL bound graphs}, this completes the proof.
\endproof

There are (at least) two natural ways to generalize the above problem. The first, considering hypergraphs, was already investigated by
Erd\H{o}s, Gy\'{a}rf\'{a}s and {\L}uczak~\cite{EGL:94}.
In another direction, while focusing on the graph case, Yuster~\cite{yuster:97} extended the definition of an $(n,k)$-graph by insisting on $s$ edges between any two parts, instead of just one edge.

For brevity, we first combine both directions in a single definition, and then discuss the respective contributions of~\cite{EGL:94} and~\cite{yuster:97}. Given a hypergraph~$G$ whose vertex set is partitioned into $V(G)=V_1\cup \dots \cup V_m$, an \defn{independent transversal} (with respect to $\Set{V_i}_{i\in [m]}$) is a set containing exactly one vertex from each~$V_i$, but not containing any edge.
Define an \defn{$(n,k,r,s)$-graph} to be an $r$-graph with $n$ disjoint parts of size $k$, such that any $r$~parts induce exactly $s$ edges, each containing one vertex from each part, and such that these $s$ edges form a matching. Let $f(k,r,s)$ be the maximum integer $n$ such that any $(n,k,r,s)$-graph has an independent transversal.
Hence, an $(n,k)$-graph is precisely an $(n,k,2,1)$-graph, and $f(k)=f(k,2,1)$.
Of course one can also dispose of the matching condition, and instead allow any $s$~edges. We want to avoid introducing more definitions for this and thus simply stick to the variant introduced by Yuster~\cite{yuster:97}. In relevant places, we will make some remarks regarding this subtlety. 

Erd\H{o}s, Gy\'{a}rf\'{a}s and {\L}uczak~\cite{EGL:94} studied the function $f(k,r,1)$, whereas Yuster~\cite{yuster:97} investigated $f(k,2,s)$ more closely.
A standard application of the Lov\'asz Local Lemma (cf.~\cite{EGL:94} and Remark~\ref{rem:LLL}) shows that
\begin{align}
	f(k,r,s)=\Omega_r\left(\sqrt[r-1]{k^r/s}\right).\label{general LLL lower bound}
\end{align}
\COMMENT{The implicit constant is $\sqrt[r-1]{\frac{(r-1)!}{\eul r}}$.}
Moreover, if we choose an $(n,k,r,s)$-graph uniformly at random, then the expected number of independent transversals is $k^n\left(1-s/k^r\right)^{\binom{n}{r}}$. For $n\gg_r \sqrt[r-1]{(k^r\log k)/s}$, this value is less than~$1$, implying that there exists an $(n,k,r,s)$-graph with no independent transversal. Hence, 
\begin{align}
	f(k,r,s)=\cO_r\left(\sqrt[r-1]{(k^r\log k)/s}\right),\label{general prob upper bound}
\end{align}
which is larger than the lower bound in~\eqref{general LLL lower bound} by a $\cO_r(\log^{1/(r-1)} k)$-factor.

As mentioned above, in the case $r=2,s=1$, Erd\H{o}s, Gy\'{a}rf\'{a}s and {\L}uczak provided an explicit construction~\cite[Proposition~6]{EGL:94}, showing that $f(k,2,1)\le (1+o(1))k^2$. However, they were unable to remove the $\log$-factor in~\eqref{general prob upper bound} for larger~$r$ (they only considered $s=1$).

Yuster~\cite[Theorem~1.3]{yuster:97} proved that $f(k,2,s)\le (1+o_k(1))k^2/s$ for a wide range of values for~$s$, namely when $\log k \ll s\ll k$, using a probabilistic construction.
Motivated by this, he conjectured that this bound holds for all values of~$s$. 
\begin{conj}[Yuster~\cite{yuster:97}]\label{conj:yuster}
$f(k,2,s)\le (1+o_k(1))\frac{k^2}{s}$ for all $1\le s\le k$.
\end{conj}

Here, we prove Conjecture~\ref{conj:yuster}, by establishing a connection to the Tur\'an problem for complete bipartite graphs. More precisely, for $n,m,r,s\in \bN$, we define the function $w(n,m;r,s)$ as the maximum integer $k$ for which there exists a bipartite graph with parts $A$ of size~$n$ and $B$ of size~$m$ such that every vertex in $A$ has degree at least~$k$, and any $r$ vertices in $A$ have at most $s$ common neighbours.

Given such a bipartite graph for suitable parameters, we can construct an $(n,k,r,s)$-graph without any independent transversal, which leads to the following result.

\begin{lemma}\label{lem:Turan reduction}
If $w(n,m;r,s)\ge k$ and $m(r-1)<n$, then $f(k,r,s)<n$.
\end{lemma}

Observe that $w(n,m;r,s)$ is a minimum degree version of the well-known \defn{Zarankiewicz function} $z(n,m;r,s)$, defined as the maximum number of edges in a bipartite graph with parts $A$ of size~$n$ and $B$ of size~$m$ such that any $r$ vertices in $A$ have less than $s$ common neighbours.
Clearly, $z(n,m;r,s+1)\ge nw(n,m;r,s)$.
Moreover, the (known) constructions attaining $z$ are usually almost regular, in which case it follows that $z(n,m;r,s+1)\approx nw(n,m;r,s)$.
The Zarankiewicz function in turn is closely related to the Tur\'an number of the complete bipartite graph~$K_{r,s}$. Determining even the order of magnitude of $\ex(n;K_{r,s})$ is a notoriously difficult problem. Fortunately, in connection with Conjecture~\ref{conj:yuster}, we have $r=2$, for which the aforementioned functions are relatively well understood.
We will prove Lemma~\ref{lem:Turan reduction} and deduce Conjecture~\ref{conj:yuster} in Section~\ref{sec:upper bounds}. 

Note that Conjecture~\ref{conj:yuster} only states an upper bound for~$f(k,2,s)$. The Lov\'asz Local Lemma (cf.~\eqref{general LLL lower bound}) yields $k^2/(2\eul s)$ as a lower bound. Yuster~\cite[Theorem~1.1]{yuster:97} improved the $(1/2\eul)$-factor to~$0.52$ in the case when $s$ is fixed and $k$ is sufficiently large. Using Theorem~\ref{thm:main general}, we can in fact show that the upper bound in Conjecture~\ref{conj:yuster} is the correct value for $f(k,2,s)$, for all~$s$. Thus, we have the following stronger form of Theorem~\ref{thm:EGL graph}.

\begin{theorem} \label{thm:yuster conj}
$f(k,2,s)= (1+o_k(1))\frac{k^2}{s}$ for all $1\le s\le k$. 
\end{theorem}

\proof
Let $G$ be any $(n,k,2,s)$-graph. Then each part has size $k$ and average degree $(n-1)s/k$, and the local degree is~$1$. Hence, if $k\ge (1+o_k(1))ns/k$, then an independent transversal exists by Theorem~\ref{thm:main general}. This implies $f(k,2,s) \ge (1-o_k(1))\frac{k^2}{s}$. Together with the confirmation of Conjecture~\ref{conj:yuster} (cf.~Section~\ref{sec:upper bounds}), this completes the proof.
\endproof

Recall that in the definition of an $(n,k,2,s)$-graph, one could also dispose of the matching condition, thus allowing any $s$~edges between two parts. Clearly, the upper bound for (the analogue of) $f(k,2,s)$ still holds. For the lower bound, note that in the proof of Theorem~\ref{thm:yuster conj}, the matching condition is only used to ensure local degree~$1$. 
Hence, for $s=o(k)$, in particular for fixed~$s$, we could still use Theorem~\ref{thm:main general} to establish the same lower bound without the matching condition, since the local degree would be obviously at most~$s$.

Theorem~\ref{thm:yuster conj} completely settles the $r=2$ case of the Erd\H{o}s--Gy\'{a}rf\'{a}s--{\L}uczak problem and its generalization due to Yuster. Perhaps the most significant open problem for $r>2$ is to get rid of the extra $\log$-factor in~\eqref{general prob upper bound}.
Unfortunately, we cannot prove that $f(k,r,1)=\cO_r(k^{r/(r-1)})$. However, when $s$ is sufficiently large compared to~$r$, then we can use known results on the Tur\'an function of $K_{r,s}$ together with Lemma~\ref{lem:Turan reduction} in order to determine the order of magnitude of~$f(k,r,s)$.

\begin{theorem}\label{thm:hypergraph upper bound}
For every fixed $r$, we have $f(k,r,s)=\Theta_r(\sqrt[r-1]{k^r/s})$ for all $(r-1)!\le s\le k$.
\end{theorem}

It seems natural to conjecture that $f(k,r,s) = \Theta_r(\sqrt[r-1]{k^r/s})$ for all~$s$.
It is perhaps less obvious what the implicit constant should be, since neither of~\eqref{general LLL lower bound} and~\eqref{general prob upper bound} give the correct answer for $r=2$.

We remark that for very large~$s$, we can obtain a more effective upper bound, using a probabilistic construction (and again Lemma~\ref{lem:Turan reduction}). The $r=2$ case of this result will also be used to prove Conjecture~\ref{conj:yuster} in this range.

\begin{theorem}\label{thm:hypergraph upper bound large s}
For every fixed $r\ge 2$ and $\eps>0$, there exists $C>0$ such that whenever $C\log k\le s\le k$, we have $f(k,r,s)\le (r-1+\eps)\sqrt[r-1]{k^r/s}$.
\end{theorem}

We will prove Theorems~\ref{thm:hypergraph upper bound} and~\ref{thm:hypergraph upper bound large s} in Section~\ref{sec:upper bounds}.

\section{Probabilistic tools}

We use some probabilistic tools in our proofs. Since they are well known, we find it more efficient to refer to them by name rather than number.

The Lov\'asz Local Lemma has already been mentioned frequently. Given some `bad' events in a probability space, it can be used to show that with positive probability, none of them happens, provided they do not depend too much on each other. To make this more precise, given events $A_1,\dots,A_n$ in a probability space, we say that a graph $\Gamma$ on $[n]$ is a \defn{dependency graph} for these events if each $A_i$ is mutually independent of all other events except those indexed by $N_\Gamma(i)\cup \Set{i}$.

\begin{theorem}[Lov\'asz Local Lemma, see~\cite{AS:08}]
Let $A_1,\dots,A_n$ be events in a probability space with dependency graph $\Gamma$. If $\prob{A_i}\le p$ for all $i\in [n]$ and $\eul p (\Delta(\Gamma)+1)\le 1$, then $\prob{\cap_{i=1}^n \overline{A_i}}>0$.
\end{theorem}

Moreover, we will need concentration results. For sums of independent variables, it will be convenient to use the well-known Chernoff--Hoeffding bound.

\begin{theorem}[Chernoff's bound, see~\cite{JLR:00}] \label{lem:chernoff}
Let $X$ be the sum of independent Bernoulli random variables (possibly with distinct expectations). Then, for all $0\le\eps \le 3/2$, $$\prob{X\neq (1\pm \eps)\expn{X}} \leq 2\exp\left(-\eps^2\expn{X}/3\right).$$
Moreover, for $t\ge 7 \expn{X}$, we have $\prob{X\ge t}\le \exp(-t)$.
\end{theorem}

For more complicated random variables, we will appeal to Talagrand's inequality. Given a product probability space $\Omega=\prod_{i=1}^n \Omega_i$ and a random variable $X\colon \Omega\to \bR$, we say that $X$ is
\begin{itemize}
\item \defn{$C$-Lipschitz} if for any $\omega\in \Omega$, changing one coordinate of $\omega$ affects the value of $X$ by at most~$C$;
\item \defn{$r$-certifiable} if for every $s$ and $\omega$ such that $X(\omega)\ge s$, there exists a set $I\In [n]$ of size at most $ rs$ such that $X(\omega')\ge s$ for every $\omega'$ that agrees with $\omega$ on the coordinates indexed by~$I$.
\end{itemize}

\begin{theorem}[Talagrand's inequality, see~\cite{MR:02}]
Suppose that $X$ is $C$-Lipschitz and $r$-certifiable. Then
$$\prob{|X-\expn{X}|>t+60C\sqrt{r\expn{X}}} \le 4\exp\left(-\frac{t^2}{8C^2r\expn{X}}\right).$$
\end{theorem}

\section{Proof of Theorem~\ref{thm:main general}}

We start by explaining the main approach, which also underpinned the proofs in~\cite{LohSudakov:07,RS:02}. 
Let $G$ be a graph with vertex partition $V(G)=V_1 \cup\dots\cup V_m$. Consider first the maximum degree setting; that is, assume that
$\Delta(G)\le \Delta$ and $|V_i|=S \ge (1+\eps)\Delta$ for all $i\in[m]$. As mentioned in the introduction, if $S \ge 2\eul \Delta$, then the
Lov\'asz Local Lemma guarantees the existence of an independent transversal. The aim is to combine this with a `nibble' algorithm. In each step, it chooses a few vertices to be included in the transversal, such that after deleting all sets $V_i$ from which a vertex has been chosen, and deleting all neighbours of chosen vertices, the ratio of remaining part sizes to maximum degree increases. After a sufficient number of such steps, this ratio will be larger than $2\eul$ so that the Lov\'asz Local Lemma is applicable to complete the transversal in a single additional step.

To see why it is plausible that this works, we sketch one such nibble step. `Activate' each part $V_i$ with some small probability~$p$, and from each activated part, select a vertex uniformly at random. Let $T$ be the set of all chosen vertices that do not have a chosen neighbour.
We add $T$ to our transversal. To avoid future conflicts, we delete all neighbours of chosen vertices. For every vertex, the probability that one of its neighbours is chosen is at most $\Delta\cdot p \cdot \frac{1}{S}\le p/(1+\eps)$ by assumption. Thus, we may expect that each part loses a $(p/(1+\eps))$-fraction of its size. On the other hand, the parts which have a vertex in $T$ can be deleted entirely. For each part, once it is activated, the probability that the selected vertex will not be in $T$ is at most $ p$ by the above. Hence, each part is deleted with probability at least $p-p^2$, so we might expect that every vertex loses at least a $(p-p^2)$-fraction of its degree. Choosing $p$ small enough, this ensures that the maximum degree shrinks by a larger factor than the size of the remaining parts, as desired.

In our proof, instead of tracking the maximum degree during the nibble process, we will track the average degree of each part, and our goal is to show that the ratio of part size to average degree increases. We observe that degrees will still shrink by roughly a $p$-fraction. (We remark that by choosing first a vertex randomly from \emph{every} part, and then discarding those of non-activated parts, we obtain a much more streamlined proof of this fact compared to~\cite{LohSudakov:07}.)
Moreover, for each part, on \emph{average}, the probability of a vertex being deleted is still at most $p/(1+\eps)$, hence we still expect each part to lose at most a $(p/(1+\eps))$-fraction of its size. To control the new average degree, we need to know in addition \emph{which} vertices are deleted from a part, since deleting predominantly low-degree vertices would increase the average degree. Fortunately, high degree vertices are more likely to have a chosen neighbour and thus to be deleted, which plays into our hands.

\subsection{Reductions}

The proof of Theorem~\ref{thm:main general} involves a number of reductions.
First of all, even though we work with average degrees, it is often needed to have a good upper bound on the maximum degree.
One simple but crucial observation is that we can assume $\Delta(G)=\cO_\eps(D)$, since an independent transversal of $G'$ in the below statement will also be an independent transversal of~$G$.

\begin{prop}\label{prop:max deg reduction}
Suppose $\eps\in (0,1)$ and let $G$ be a graph with vertex partition $V(G)=V_1 \cup\dots\cup V_m$. Suppose that $|V_i|\ge (1+\eps)D$ and $\bar{d}_G(V_i)\le D$ for each $i\in[m]$. Then there exist $D'\ge D$ and subsets $V_i'\In V_i$ of size at least $ (1+\eps/2)D'$ such that, letting $G':=G[V_1',\dots,V_m']$, we have $\bar{d}_{G'}(V_i')\le D'$ for each $i\in[m]$, and $\Delta(G')\le 8D'/\eps$. 
\end{prop}

\proof
For each $i\in [m]$, let $B_i$ be the set of $v\in V_i$ with $d_G(v)>8D/\eps$. Then
\begin{align*}
|B_i| \cdot 8D/\eps \le \sum_{v\in V_i}d_G(v)\le D|V_i|,
\end{align*}
whence $|B_i|\le \eps|V_i|/8$. Define $V_i':=V_i\sm B_i$ and $G':=G[V_1',\dots, V_m']$.
Then $$|V_i'|\ge (1-\eps/8)|V_i|\ge (1-\eps/8)(1+\eps)D \mbox{ and }\bar{d}_{G'}(V_i')\le \frac{|V_i|}{|V_i'|}\bar{d}_G(V_i)\le D/(1-\eps/8).$$
Let $D':=D/(1-\eps/8)$.
Then we have $\bar{d}_{G'}(V_i')\le D'$ and $|V_i'|\ge (1-\eps/8)^2(1+\eps)D'\ge (1+\eps/2)D'$ for all $i\in[m]$. Moreover, $\Delta(G')\le 8D/\eps \le 8D'/\eps$.
\endproof

A further reduction concerns the local degree. Recall that in the nibble step sketched above, each part is deleted with probability roughly~$p$, and we want to conclude that each vertex loses this fraction of its degree. To apply concentration inequalities, each part must only have a small effect on the total degree. In other words, we need small local degree. We will analyze the nibble procedure to prove the following result.

\begin{theorem}\label{thm:main}
For all $\eps>0$ the following holds for sufficiently large~$D$.
Let $G$ be a graph with vertex partition $V(G)=V_1 \cup\dots\cup V_m$ and local degree at most $ \log^2 D$. Suppose that $|V_i|\ge (1+\eps)D$ and $\bar{d}_G(V_i)\le D$ for each $i\in[m]$. Then there exists an independent transversal.
\end{theorem}

Theorem~\ref{thm:main} would already be sufficient for our purposes here, as all our applications have local degree~$1$. 
The general case when the local degree is $o(D)$ can be reduced to Theorem~\ref{thm:main}, by finding a suitable induced subgraph $G'\In G$ which meets the conditions of Theorem~\ref{thm:main}.
First, observe that if the local degree is only $D^{1-\gamma}$ for some small constant $\gamma$, then there is a simple one-step reduction. Retain every vertex with probability~$D^{-1+\gamma}$. Then each part has expected residual size at least $(1+\eps)D^{\gamma}$ and average degree at most~$D^{\gamma}$. Moreover, for each vertex, we only expect at most one neighbour in every part. Since dependencies are polynomial in~$D$, we can use Chernoff's bound and the Lov\'asz Local Lemma to show that with positive probability, all parts have average degree at most $ (1+\eps/4)D^{\gamma}=:D'$ and size at least $ (1+\eps/2)D'$, and the local degree is at most $ \log D \log\log D \le \log^2 D'$.

If the local degree is $\gamma D$ as assumed in Theorem~\ref{thm:main general}, then this one-step reduction fails. However, a more careful sequential procedure still works.
This was also done in the maximum degree setting in~\cite{LohSudakov:07} (even reducing to constant local degree), with a method inspired by~\cite{alon:92}. Since a very similar proof also works here, we omit the details (they may be found in~\cite[appendix]{KK:ta}).

\subsection{Setting the stage}

We now prove Theorem~\ref{thm:main}. By Proposition~\ref{prop:max deg reduction}, we can assume that $\Delta:=\Delta(G)= \cO(D/\eps)$. Moreover, by repeatedly deleting vertices of largest degree, we can also assume that $|V_i|=\lceil(1+\eps)D\rceil$ for all $i\in [m]$, and we may clearly assume that each part $V_i$ is independent.
We also assume that $\eps$ is sufficiently small and $D$ is sufficiently large.

Set $p:=1/\log^3 D$.
Initiate $S(0)=(1+\eps)D$ and $D(0)=D$ and recursively define
\begin{align}
S(t+1)=\left(1-\frac{p}{1+3\eps/4}\right)S(t), \qquad D(t+1)=\left(1-\frac{p}{1+\eps/4}\right)D(t).\label{def trajectories}
\end{align}

We will inductively show for all integers $0\le t\le t^*:=\lceil \frac{10}{\eps p}\rceil$ that there exist
\begin{enumerate}[label=\rm{(P\arabic*)}]
\item a set of indices $I(t)\In [m]$;\label{ind ind}
\item an independent transversal $T(t)$ of $G[\cup_{i\in [m]\sm I(t)}V_i ]$;\label{ind trans}
\item a subset $V_i(t)\In V_i \sm N(T(t))$ for each $i\in I(t)$ such that $|V_i(t)|\ge S(t)$  and $\bar{d}_{G(t)}(V_i(t)) \le D(t)$, where $G(t):=G[\cup_{i\in I(t)}V_i(t)]$.\label{ind sets}
\end{enumerate}

Here, we write $N(v)$ for the neighbourhood of $v$ in $G$ and $N(S):=\cup_{v\in S}N(v)$.

Note that at any stage, since the subsets $V_i(t)$ do not contain any neighbours of $T(t)$, it is sufficient to find an independent transversal of $G(t)$ in order to complete $T(t)$ to an independent transversal of~$G$.
Fortunately, in $G(t^*)$, the remaining parts have size at least $ S(t^*)$ and average degree at most $ D(t^*)$. A routine calculation reveals that $S(t^*)\ge 2\eul D(t^*)$, with room to spare. Indeed,
\begin{align*}
\frac{D(t^*)}{S(t^*)} &\le \left(\frac{1-\frac{p}{1+\eps/4}}{1-\frac{p}{1+3\eps/4}}\right)^{t^*} \le \left(1-\frac{p\eps}{5}\right)^{\frac{10}{\eps p}} \le \eul^{-2} \le \frac{1}{2\eul}.
\end{align*}
\COMMENT{$\frac{1-\frac{p}{(1+\eps/4)}}{1-\frac{p}{(1+3\eps/4)}} = 1-\frac{p(\frac{1}{(1+\eps/4)}-\frac{1}{(1+3\eps/4)})}{1-\frac{p}{(1+3\eps/4)}} \le 1-p(\frac{1}{(1+\eps/4)}-\frac{1}{(1+3\eps/4)})$ and $\frac{1}{(1+\eps/4)}-\frac{1}{(1+3\eps/4)}=\frac{\eps}{2(1+\eps/4)(1+3\eps/4)}\ge \eps/5$}
Thus, we can complete the transversal by applying the following to $G(t^*)$.
\begin{lemma} \label{lem:LLL app}
Let $G$ be a graph with vertex partition $V(G)=V_1 \cup\dots\cup V_m$. Suppose that $|V_i|\ge 2\eul D$ and $\bar{d}_G(V_i)\le D$ for each $i\in[m]$. Then there exists an independent transversal.
\end{lemma}

\proof
We may clearly assume that each $V_i$ is an independent set. Moreover, by repeatedly deleting vertices of largest degree within a part, we can assume that $|V_i|=\lceil 2\eul D\rceil=:n$ for all $i\in [m]$.
For each $i\in[m]$, pick a vertex $v_i\in V_i$ uniformly at random, independently from all the other choices. Let $T:=\Set{v_1,\dots,v_m}$. Our goal is to show that $T$ is independent with positive probability. For an edge $e\in E(G)$, let $A_e$ be the event that both of its endpoints are in~$T$. Clearly, $\prob{A_e}= n^{-2}$. We define a graph $\Gamma$ on $E(G)$ by joining $e$ and $e'$ if one endpoint of $e$ and one endpoint of $e'$ lie in the same part. Observe that $\Gamma$ is a dependency graph for the events $\Set{A_e}_{e\in E(G)}$. Fix an edge $e$ with endpoints in $V_i$ and~$V_j$. Then
$$d_\Gamma(e) \le |V_i|\bar{d}_G(V_i) + |V_j|\bar{d}_G(V_j) - 2. $$ Hence, $\Delta(\Gamma) \le 2nD-2$. Since $e\cdot n^{-2}\cdot 2nD\le 1$, the Lov\'asz Local Lemma applies and we can infer that with positive probability, none of the events $A_e$ happens.
\endproof

\begin{remark}\label{rem:LLL}
The above result is essentially due to Alon~\cite{alon:88}. He proved it with $25$ in place of $2\eul$, but see \cite[Chapter~5]{AS:08}. Moreover, both references state the result in terms of the maximum degree, but the same proof works in the average degree setting. The proof works as well for hypergraphs and can be used to deduce \eqref{general LLL lower bound}. We also note that the $s=1$ case of \eqref{general LLL lower bound} is already contained in~\cite{EGL:94}.
\end{remark}

Note also that \ref{ind ind}--\ref{ind sets} hold for $t=0$, with $I(0)=[m]$, $T(0)=\emptyset$ and $V_i(0)=V_i$ for all $i\in [m]$.
It remains to show that if \ref{ind ind}--\ref{ind sets} hold for $t<t^*$, then they also hold for~$t+1$.
Note that for the entire range of~$t$, we have $S(t)=\Theta_\eps(D)$ and $D(t)=\Theta_\eps(D)$. We will use this fact throughout the rest of the proof.

Assume that $I(t)$, $T(t)$, $\Set{V_i(t)}_{i\in I(t)}$, $G(t)$ are given.
We proceed as follows. For each $i\in I(t)$, we activate $i$ with probability~$p$. Let $J\In I(t)$ be the random subset of activated parts. For each $i\in J$, we pick a vertex $v_i\in V_i(t)$ uniformly at random, and let $T':=\set{v_i}{i\in J}$ be the set of all selected vertices. Note that $T'$ might induce edges.
Let $\hat{J}$ be the set of $i\in J$ for which $v_i$ does not have a neighbour in $T'$, and let $\hat{T} :=\set{v_i}{i\in \hat{J}}$. Clearly, $\hat{T}$ is independent.
Let $T(t+1):=T(t)\cup \hat{T}$ and $I(t+1):=I(t)\sm \hat{J}$.
Then \ref{ind ind} and \ref{ind trans} hold.

We need to delete all neighbours of vertices in $\hat{T}$ from the graph. In fact, to obtain better control on certain variables, we potentially delete more vertices than that. 
First note that for every vertex $v\in V(G(t))$, each of its remaining neighbours is included in $T'$ with probability at most $ p/S(t)$, whence 
\begin{align*}
	\prob{|N(v) \cap T'|\ge 1} \le d_{G(t)}(v)\frac{p}{S(t)} = \cO_\eps(p). 
\end{align*}
This implies
\begin{align*}
\prob{N(v) \cap T'=\emptyset}   \ge 1-d_{G(t)}(v)\frac{p}{S(t)} =: p_v.
\end{align*}
We define an artificial Bernoulli random variable $B_v$, independent of all other variables, such that
\begin{align}
\prob{N(v) \cap T'=\emptyset \wedge B_v=1} = p_v. \label{truncate}
\end{align}

Now, for each $i\in I(t)$, define $V_i(t+1)$ as the set of all $v\in V_i(t)$ with $N(v) \cap T'=\emptyset$ and $B_v=1$.
Since $\hat{T}\In T'$, this ensures that $V_i(t+1)\In V_i \sm N(T(t+1))$ for each $i\in I(t+1)$. 
It remains to establish the conditions on the size and average degree of~$V_i(t+1)$ in~\ref{ind sets}.

\subsection{Two key claims}

For this, we utilize two claims. Let $C:=\log^2 D$ be the upper bound on the local degree of~$G$. The proofs below use concentration inequalities and are the only places where the assumption on the local degree is needed.

The first claim asserts that given any subset of one part, with high probability, the number of vertices that remain is concentrated around its expectation. Note that the `error term' can be larger than $|U|$ if $U$ is small, making the statement trivial. However, it will be convenient later on not to assume a lower bound on the size of~$U$.

\begin{claim} \label{claim:set sizes}
For any index $i\in I(t)$ and subset $U\In V_i(t)$, we have $$\prob{|U\cap V_i(t+1)|\neq \sum_{v\in U}p_v \pm C\log D\sqrt{p|U|}} \le \exp(-\log D\log\log D).$$
\end{claim}

\claimproof
From \eqref{truncate}, we immediately have $\expn{|U\cap V_i(t+1)|}=\sum_{v\in U}p_v$. We now establish concentration. Let $R=|U\sm V_i(t+1)|$ be the number of vertices removed from $U$. Observe that $R$ is $1$-certifiable. (This is the reason why we consider $R$ instead of $|U\cap V_i(t+1)|$.) Indeed, for any vertex $v$ removed, there must either be a neighbour of $v$ in $T'$, or we have $B_v=0$, and these events certify the removal of~$v$.
Moreover, $R$ is clearly $C$-Lipschitz due to the bound on the local degree.
Note that $\expn{R}=\sum_{v\in U}(1-p_v) = \cO_\eps(p|U|)$. Thus, applying Talagrand's inequality with $t=C\log D\sqrt{p|U|}/2$, we obtain
\begin{align*}
\prob{|R-\expn{R}|>2t}\le \exp(-\log D\log\log D),
\end{align*}
since $60C\sqrt{\expn{R}}=\cO_\eps\left(C\sqrt{p|U|}\right)\le t$ and $t^2/(C^2\expn{R})=\Omega_\eps(\log^2 D)$.
This proves the claim as $|U\cap V_i(t+1)|-\sum_{v\in U}p_v=\expn{R}-R$.
\endclaimproof

The second claim asserts that each vertex loses the right proportion of its neighbours. This corresponds to (\rm{\bf ii}) in~\cite{LohSudakov:07}. We give a shorter proof here.

\begin{claim} \label{claim:degrees}
For any vertex $v\in V(G(t))$ with $d_{G(t)}(v)\ge \log^{20} D$, we have $$\prob{d_{G(t+1)}(v) \ge  \left(1-\frac{p}{1+\eps/5}\right) d_{G(t)}(v)} \le \exp(-\log D\log\log D).$$
\end{claim}

\claimproof
Fix a vertex $v\in V(G(t))$ with $d_{G(t)}(v)\ge \log^{20} D$.
For $i\in I(t)$, let $a_i:=|N(v)\cap V_i(t)|$, and let $K$ be the set of $i\in I(t)$ for which $a_i>0$. 
Define the random variable $Z=\sum_{i\in K}a_i\IND_{i\in \hat{J}}$ which counts the number of neighbours $v$ loses because their entire part is deleted. Observe that $$d_{G(t+1)}(v)\le d_{G(t)}(v) -Z = \sum_{i\in K}a_i\IND_{i\notin \hat{J}} =:\bar{Z}.$$ Hence, in order to prove the claim, it suffices to show that
\begin{align}
	\prob{\bar{Z} \ge  \left(1-\frac{p}{1+\eps/5}\right) d_{G(t)}(v)} \le \exp(-\log D\log\log D).\label{neighbour deletion}
\end{align}
We now describe an alternative way of sampling~$T'$. In the first phase, choose a vertex $v_i$ uniformly at random from $V_i(t)$ for \emph{every} $i\in I(t)$. In the second phase, activate each $i\in I(t)$ with probability~$p$. Then $T'$ is exactly the set of $v_i$ for which $i$ has been activated. The advantage of this two-round exposure is that, typically, the chosen vertices $v_i$ induce a graph with very small maximum degree. Conditioning on this when we activate the parts, we have a much better Lipschitz condition which allows us to apply Talagrand's inequality.

Now, first randomly choose the vertices~$v_i$. Let $\bar{K}$ be the set of $i\in I(t)$ for which $i\in K$ or there is $j\in K$ such that $e_G(V_i,V_j)>0$. Note that $|\bar{K}|=\cO_\eps(D^3)$ and thus $|\cup_{i\in \bar{K}}V_i(t)|=\cO_\eps(D^4)$. 
Let $\cE$ be the event that every vertex in $\cup_{i\in \bar{K}}V_i(t)$ has at most $\log^2 D$ neighbours in $\set{v_i}{i\in I(t)}=:T''$.
We first claim that 
\begin{align}
	\prob{\overline{\cE}}\le \frac{1}{2}\exp(-\log D\log\log D).\label{pre conditioning}
\end{align}
Indeed, for any vertex $w$ in $\cup_{i\in \bar{K}}V_i(t)$ and $j\in I(t)$, the probability that $v_j\in N(w)$ is $$|N(w)\cap V_j(t)|/|V_j(t)|\le |N(w)\cap V_j(t)|/S(t).$$ Hence, the expected number of $j$ for which $v_j\in N(w)$ is at most $d_{G(t)}(w)/S(t)=\cO_\eps(1)$.
By Chernoff's bound, the probability that $|N(w)\cap T''|> \log^2 D$ is at most $\exp(-\log^2 D)$. A union bound over the $\cO_\eps(D^4)$ choices for $w$ implies~\eqref{pre conditioning}.

Next, we show that $\bar{Z}$ is unlikely to be too large if $\cE$ holds. For this, fix any choice $T''$ of $\set{v_i}{i\in I(t)}$ for which $\cE$ holds. For $i\in I(t)$, let $Q_i$ be the set of $j\in I(t)$ with $v_iv_j\in E(G)$. By assumption, we have $|Q_i|\le \log^2 D$ for all $i\in \bar{K}$.
Now, we randomly activate the parts. Observe that $i\in \hat{J}$ if and only if $i$ is activated but no $j\in Q_i$ is activated. Hence, for $i\in K$, 
$$p\ge \cprob{i\in \hat{J}}{T''} \ge p(1-p|Q_i|)\ge p(1-1/\log D),$$ where the last inequality holds by our choice of~$p$.
Thus, $\cexpn{\bar{Z}}{T''} \le (1-p+p/\log D)d_{G(t)}(v)$. Moreover, $\bar{Z}$ only depends on the activation of $i$ if $i\in \bar{K}$, and whether such $i$ is activated or not only affects the events $\Set{j\in \hat{J}}$ with $j\in Q_i\cup \Set{i}$. Since $a_j\le C$ for all $j$, we infer that $\bar{Z}$ is Lipschitz with constant $C(1+|Q_i|)\le 2\log^4 D$.
Finally, observe that $\bar{Z}$ is $1$-certifiable, because the event $\Set{i\notin \hat{J}}$ is certified by either the non-activation of $i$ or the activation of some $j\in Q_i$.
We now apply Talagrand's inequality with $t=pd_{G(t)}(v)/\log D$. Noting that $60C\sqrt{\cexpn{\bar{Z}}{T''}}\le 60C\sqrt{d_{G(t)}(v)} \le t$ and $t^2/(C^2\cexpn{\bar{Z}}{T''})=\Omega(\log^2 D)$, where we use the assumption $d_{G(t)}(v)\ge \log^{20} D$, we have 
$$\cprob{\bar{Z} \ge  (1-p+3p/\log D))d_{G(t)}(v)}{T''} \le \frac{1}{2}\exp(-\log D\log\log D).$$ Since this bound holds for any choice $T''$ that satisfies~$\cE$, we conclude that
$$\cprob{\bar{Z} \ge  \left(1-\frac{p}{1+\eps/5}\right) d_{G(t)}(v)}{\cE}\le \frac{1}{2}\exp(-\log D\log\log D).$$
Adding this inequality to~\eqref{pre conditioning} implies \eqref{neighbour deletion} and thus the claim.
\endclaimproof

\subsection{Putting everything together}

We now use Claims~\ref{claim:set sizes} and~\ref{claim:degrees} to complete the proof. Note that it is enough to show that for each $i\in I(t)$, the `bad event' $A_i$ that one of $|V_i(t+1)| < S(t+1)$ or $\bar{d}_{G(t+1)}(V_i(t+1))>D(t+1)$ happens, has probability at most $\cO(D)\exp(-\log D\log\log D)$. Indeed, although we cannot simply use a union bound since we do not assume any bound on the number of parts~$m$, we can easily use the Lov\'asz Local Lemma to show that with positive probability, no $A_i$ happens. Let $\Gamma$ be the graph on $I(t)$ where $i,j$ are joined if there is at least one edge between $V_i,V_j$. Observe that the event $A_i$ is fully determined by the random choices involving parts $V_j$ with $\mathrm{dist}_\Gamma(i,j)\le 2$.
Thus, $\Gamma^4$, the graph on $I(t)$ where $i,j$ are joined if their distance in $\Gamma$ is at most~$4$, is a dependency graph for the events $\Set{A_i}_{i\in I(t)}$, since if $ij\notin E(\Gamma^4)$, then $A_i,A_j$ are determined by disjoint sets of random choices.
Moreover, $\Delta(\Gamma)=\cO_\eps(D^2)$, implying $\Delta(\Gamma^4)=\cO_\eps(D^8)$. Hence, the condition of the Lov\'asz Local Lemma is fulfilled and so there is an outcome which satisfies~\ref{ind sets}.

Fix $i\in I(t)$. The desired bound for $|V_i(t+1)|$ will follow relatively easily from Claim~\ref{claim:set sizes}. 
In order to estimate $\bar{d}_{G(t+1)}(V_i(t+1))$, we proceed in two steps. Note that using Claim~\ref{claim:degrees}, we know that every vertex loses the right proportion of its degree. Thus, the average of $d_{G(t+1)}(v)$ over all $v\in V_i(t)$ will indeed be bounded by~$D(t+1)$. However, this alone is not enough, since deleting vertices from $V_i(t)$ could increase the average degree if we were to delete predominantly low-degree vertices. But of course, high-degree vertices are more likely to be deleted, which is exactly why the proof works. 

First, we group the vertices of $V_i(t)$ according to their degree. Set $b:=\lceil D^{2/3}\rceil$ and for $j\in \bN\cup \Set{0}$, let
$$X_j:=\set{v\in V_i(t)}{jb\le d_{G(t)}(v) <(j+1)b }$$
and $X_j':=X_j \cap V_i(t+1)$. Note that $X_j\neq \emptyset$ only if $j=\cO(\Delta/b)=\cO_\eps(D^{1/3})$. 
Clearly, $|V_i(t)|=\sum_j |X_j|$ and $|V_i(t+1)|=\sum_j |X_j'|$.

Using Claims~\ref{claim:set sizes} and~\ref{claim:degrees} and a union bound, we infer that the following hold with probability at least $1-\cO(D)\exp(-\log D\log\log D)$:
\begin{enumerate}[label=\rm{(\arabic*)}]
\item $|V_i(t+1)| \ge \sum_{v\in V_i(t)}p_v - \log^2 D \sqrt{|V_i(t)|}$;\label{induction global size}
\item $|X_j'| = \sum_{v\in X_j}p_v \pm \log^2 D\sqrt{|X_j|}$ for all~$j$;\label{induction local size}
\item $d_{G(t+1)}(v) \le  \left(1-\frac{p}{1+\eps/5}\right) d_{G(t)}(v)$ for all $v\in V_i(t)\sm X_0$.\label{induction degrees}
\end{enumerate}
It thus only remains to show that assuming the above, we can deduce $|V_i(t+1)| \ge S(t+1)$ and $\bar{d}_{G(t+1)}(V_i(t+1))\le D(t+1)$.

Note that 
\begin{align*}
\sum_{v\in V_i(t)}d_{G(t)}(v) =|V_i(t)|\bar{d}_{G(t)}(V_i(t)) \le |V_i(t)|D(t)\le \frac{|V_i(t)| S(t)}{1+\eps},
\end{align*}
where we have used that $D(t)/S(t) \overset{\eqref{def trajectories}}{\le} D(0)/S(0) = 1/(1+\eps)$. Hence,
\begin{align*}
\sum_{v\in V_i(t)} p_v &= |V_i(t)| - \frac{p}{S(t)}\sum_{v\in V_i(t)}d_{G(t)}(v) \ge \left(1-\frac{p}{1+\eps}\right)|V_i(t)|.
\end{align*}

Therefore, \ref{induction global size} immediately implies
\begin{align*}
|V_i(t+1)| \ge \left(1-\frac{p}{1+\eps}-\frac{\log^2 D}{\sqrt{|V_i(t)|}}\right)|V_i(t)|\ge \left(1-\frac{p}{1+3\eps/4}\right)|V_i(t)|
\end{align*}
and thus $|V_i(t+1)|\ge S(t+1)$.

We now turn to the average degree. As mentioned above, we proceed in two steps. First, we investigate how the deletion of vertices in $V_i(t)$ affects the average degree.
Our goal is to show that 
\begin{align}
	\frac{\sum_{v\in V_i(t+1)}d_{G(t)}(v)}{|V_i(t+1)|} \le \left(1+D^{-1/5}\right)D(t).\label{crucial degree sum}
\end{align}

Define $q:=pb/S(t)$ and note that we have $jq\le 1/2$ for all~$j=\cO_\eps(D^{1/3})$. In particular, we have $\sum_{j}\left(1-jq\right)|X_j|=\Theta_\eps(D)$.
For $v\in X_j$, we have $p_v = 1-(j\pm 1)b\frac{p}{S(t)}= 1-jq \pm D^{-1/3}$. Thus, together with~\ref{induction local size}, $$|X_j'|= \left(1-jq\right)|X_j| \pm \left(|X_j|/D^{1/3}+\log^2 D\sqrt{|X_j|}\right).$$
By concavity of the square root function, the additive errors accumulate to at most 
$$\sum_{j}\left(|X_j|/D^{1/3}+\log^2 D\sqrt{|X_j|}\right) = \cO_\eps\left(D^{2/3}\log^2 D\right).$$
(Recall that the number of $j$ with $X_j\neq \emptyset$ is $\cO_\eps(D^{1/3})$.)
We deduce that
\begin{align}
  |V_i(t+1)| &= \sum_{j}|X_j'|  \ge \left(1-D^{-1/4}\right)\sum_{j}\left(1-jq\right)|X_j|   \label{vertex survive}
\end{align}
	and
\begin{align}
	\sum_{v\in V_{i+1}(t)} d_{G(t)}(v) &\le \sum_{j}|X_j'|(j+1)b \le \sum_{j}\left(1-jq\right)|X_j|jb + D^{7/4}. \label{degree survive}
\end{align}

The following is the key estimate, where the left hand side is an approximation of the average degree of $V_i(t+1)$ and the right hand side is an approximation of the average degree of~$V_i(t)$.
\begin{align}
	\bar{d}':=\frac{\sum_{j}\left(1-jq\right)|X_j|jb}{\sum_{j}\left(1-jq\right)|X_j|} \le \frac{\sum_{j}|X_j|jb}{\sum_{j}|X_j|} =:\bar{d} \label{key average estimate}
\end{align}
In order to prove~\eqref{key average estimate}, consider the equivalent inequality
$$\sum_{j,k}|X_j|\left(1-kq\right)|X_k|kb \le \sum_{j,k}|X_j|jb\left(1-kq\right)|X_k|.$$
Now double each side by adding the same terms, with the roles of $j,k$ swapped. That the inequality thus obtained holds is easily seen using that $(j-k)^2\ge 0$ for all~$j,k$.

Combining  \eqref{vertex survive} and \eqref{degree survive}, we obtain
\begin{align*}
	\frac{\sum_{v\in V_i(t+1)}d_{G(t)}(v)}{|V_i(t+1)|} &\le \frac{\sum_{j}\left(1-jq\right)|X_j|jb + D^{7/4}}{(1-D^{-1/4})\sum_{j}\left(1-jq\right)|X_j|} \le \left(1+2D^{-1/4}\right)\bar{d}'  + \cO_\eps(D^{3/4}).
\end{align*}
Moreover, we clearly have
\begin{align*}
	\bar{d}  &\le \frac{\sum_{v\in V_i(t)}d_{G(t)}(v)}{|V_i(t)|} \le D(t).
\end{align*}
In order to establish~\eqref{crucial degree sum}, we combine the above with~\eqref{key average estimate}, yielding
\begin{align*}
	\frac{\sum_{v\in V_i(t+1)}d_{G(t)}(v)}{|V_i(t+1)|} &\le \left(1+2D^{-1/4}\right)\bar{d}'  + \cO_\eps\left(D^{3/4}\right) \\
	&\le \left(1+2D^{-1/4}\right) D(t)  + \cO_\eps\left(D^{3/4}\right) \\
	&\le \left(1+D^{-1/5}\right)D(t),
\end{align*}
as claimed.

Finally, using \ref{induction degrees} and \eqref{crucial degree sum} and the trivial bound $|X_0'|\le |V_i(t)|$, we can verify that
\begin{align*}
\sum_{v\in V_i(t+1)} d_{G(t+1)}(v) &\le |X_0'|b + \sum_{v\in V_i(t+1)\sm X_0} d_{G(t+1)}(v) \\
      &\le \cO\left(D^{5/3}\right) + \left(1-\frac{p}{1+\eps/5}\right) \sum_{v\in V_i(t+1)\sm X_0} d_{G(t)}(v) \\
			&\le \cO\left(D^{5/3}\right) + \left(1-\frac{p}{1+\eps/5}\right)\left(1+D^{-1/5}\right)|V_i(t+1)|D(t) \\
			&\le \left(1-\frac{p}{1+\eps/5} + \cO\left(D^{-1/3}\right) + D^{-1/5} \right)|V_i(t+1)|D(t).
\end{align*}
Dividing by $|V_i(t+1)|$ and absorbing lower order terms, we can conclude that $\bar{d}_{G(t+1)}(V_i(t+1)) \le \left(1-\frac{p}{1+\eps/4}\right)D(t)=D(t+1)$, thus completing the proof of Theorem~\ref{thm:main}.

\section{Upper bound constructions}\label{sec:upper bounds}

We now prove Conjecture~\ref{conj:yuster} and Theorems~\ref{thm:hypergraph upper bound} and~\ref{thm:hypergraph upper bound large s}. To this end, we first need Lemma~\ref{lem:Turan reduction}, which establishes a connection between $f(k,r,s)$ and the Tur\'an problem for complete bipartite graphs. 

\lateproof{Lemma~\ref{lem:Turan reduction}}
Suppose that $w(n,m;r,s)\ge k$ and $m(r-1)<n$. Hence, there exists a bipartite graph $H$ with parts $A$ of size~$n$ and $B$ of size~$m$ such that every vertex in $A$ has degree at least~$k$, and any $r$ vertices in $A$ have at most $s$ common neighbours.
Our aim is to construct an $(n,k,r,s)$-graph with no independent transversal.

We define the \defn{neighbourhood incidence $r$-graph} $G$ induced by $H$ as follows: For every vertex $a\in A$, we define the vertex set $V_a:=\Set{a}\times N_H(a)$. The vertex set of $G$ is the disjoint union of the sets~$\Set{V_a}_{a\in A}$. For $r$ distinct vertices $a_1,\dots,a_r\in A$, we add the edge $\Set{(a_1,b_1),\dots,(a_r,b_r)}$ whenever $b_1=b_2=\dots=b_r=:b$. Note that this implies that $b$ is a common neighbour of $a_1,\dots,a_r$ in~$H$. Hence, the number of edges in $G$ induced by the parts $V_{a_1},\dots,V_{a_r}$ is $|\cap_{i=1}^r N_H(a_i)|\le s$. Moreover, these edges clearly form a matching.

Suppose, for the sake of contradiction, that $G$ has an independent transversal. This means that for each $a\in A$, there is a unique $b_a\in B$ such that $(a,b_a)$ is contained in the transversal. Moreover, for the transversal to be independent, each $b\in B$ can be used as $b_a$ for at most $r-1$ vertices $a\in A$. Hence, $|B|(r-1)\ge |A|$, contradicting our assumption that $m(r-1)<n$. Therefore, $G$ cannot have an independent transversal.

Note that $G$ is not necessarily an $(n,k,r,s)$-graph. However, by deleting vertices to make each part have size exactly~$k$, and adding some more edges such that any $r$ parts induce a matching of size~$s$, we can easily find an $(n,k,r,s)$-graph which still has no independent transversal. Thus, $f(k,r,s)<n$.
\endproof

In order to use Lemma~\ref{lem:Turan reduction}, we need good lower bounds for $w(n,m;r,s)$. For small $s$, we use known constructions of $K_{r,s}$-free graphs, which are mostly explicit and geometric in nature. 
For large~$s$, we use a simple probabilistic construction which yields the following lemma.

\begin{lemma}\label{lem:prob Turan bound}
For all $\eps>0$ and $r\ge 2$ there is $C>0$ such that whenever $m\ge s\ge C \log n$ and $n\ge r$, we have $w(n,m;r,s)\ge (1-\eps)s^{1/r}m^{1-1/r}$.
\end{lemma}

\proof
This is a standard application of the probabilistic method. Take vertex sets $A$ of size $n$ and $B$ of size $m$, and randomly include every edge between $A$ and $B$ independently with probability $p:=(1-\eps/2)(s/m)^{1/r}$. (We can clearly assume that $\eps$ is sufficiently small. In particular, we assume that $\eps r\le 1$.)

For a subset $S\In A$, let $X_S$ be the random variable counting the common neighbours of the vertices in~$S$. Clearly, $X_S$ is the sum of $m$ independent Bernoulli random variables, each with expectation $p^{|S|}$. Hence, Chernoff's bound implies that $X_S=(1\pm \eps^2)p^{|S|} m$ with probability at least $1-2\exp(-\eps^4 p^{|S|} m/3)$. For $|S|\le r$, note that $$p^{|S|} m\ge p^r m\ge (1-\eps r/2)s \ge \frac{1}{2}C \log n.$$ For sufficiently large $C$, say $C=10r\eps^{-4}$, a union bound implies that with positive probability, we have $X_S=(1\pm \eps^2)p^{|S|} m$ simultaneously for all $S\In A$ of size at most~$r$.

Thus, such an outcome exists. It is then easy to check that any $r$ vertices in $A$ have at most $(1+\eps^2)p^r m\le s$ common neighbours, as required, and each vertex in $A$ has degree at least $(1-\eps^2)pm\ge (1-\eps)s^{1/r}m^{1-1/r}$, hence establishing this as a lower bound for $w(n,m;r,s)$.
\endproof

Theorem~\ref{thm:hypergraph upper bound large s} easily follows from Lemmas~\ref{lem:Turan reduction} and~\ref{lem:prob Turan bound}.

\lateproof{Theorem~\ref{thm:hypergraph upper bound large s}}
Let $n:=\lfloor (r-1+\eps)(k^r /s)^{1/(r-1)}\rfloor$ and $m:=\lfloor n/(r-1)\rfloor-1$. Choosing $C$ large enough (note that $\log n=\Theta_r(\log k)$),  Lemma~\ref{lem:prob Turan bound} implies that $w(n,m;r,s)\ge (1-\eps^2)s^{1/r}m^{1-1/r} \ge k$. With Lemma~\ref{lem:Turan reduction}, we obtain $f(k,r,s)<n$.
\endproof

We now state the results on $K_{r,s}$-free graphs and the related Zarankiewicz problem that we need. In terms of upper bounds, the famous K\H{o}v\'ari--S\'os--Tur\'an theorem~\cite{KST:54} states that
\begin{align}
	2\ex(n,K_{r,s}) \le z(n,n;r,s) \le (1+o(1))(s-1)^{1/r}n^{2-1/r}
\end{align}
for fixed $r\le s$.
It is widely believed that the order of magnitude in this upper bound is optimal. A major challenge is to devise constructions which establish the matching lower bounds.
First, consider $r=2$. For a long time, the only known asymptotic result was that $\lim_{n\to \infty}\ex(n,K_{2,2})n^{-3/2}=\frac{1}{2}$, proved simultaneously and independently by Brown~\cite{brown:66} and Erd\H{o}s, R\'{e}nyi and S\'{o}s~\cite{ERS:66}.
F\"{u}redi~\cite{furedi:96} vastly extended this by showing that $\lim_{n\to \infty}\ex(n,K_{2,s})n^{-3/2}=\frac{1}{2}\sqrt{s-1}$ for any fixed~$s$.
For Zarankiewicz's problem, the corresponding result had already been shown earlier by M\"ors. The following theorem, which we will use to prove Conjecture~\ref{conj:yuster}, implies that $$\lim_{n\to \infty}z(n,n,2,s)n^{-3/2}=(s-1)^{1/2}$$ for any fixed~$s$.

\begin{theorem}[M\"ors~\cite{mors:81}]\label{thm:mors}
For $s\in \bN$ and a prime $p\equiv 1\mod{s}$, there exists a bipartite graph with parts $A,B$ each of size $p(p-1)/s$, such that every vertex in $A$ has degree $p-1$ and every two vertices in $A$ have at most $s$ common neighbours.
\end{theorem}

For $r>2$, much less is known. Brown~\cite{brown:66} showed that $\ex(n,K_{3,3})=\Theta(n^{5/3})$, but the order of magnitude of $\ex(n,K_{r,r})$ is still unknown for every $r\ge 4$.
In a breakthrough, Koll\'{a}r, R\'{o}nyai and Szab\'{o}~\cite{KRS:96} introduced `norm-graphs' to prove that $\ex(n,K_{r,s})=\Theta_{r,s}(n^{2-1/r})$ when $s>r!$. 
Their construction was slightly improved by Alon, R\'{o}nyai and Szab\'{o}~\cite{ARS:99}, thus showing that the result already holds when $s>(r-1)!$. Moreover, their modification reveals the right dependence on~$s$: for fixed $r$ and $s>(r-1)!$, it holds that $$\ex(n,K_{r,s})=\Theta_r((s-1)^{1/r}n^{2-1/r}).$$ This follows from the next result, which we will use to prove Theorem~\ref{thm:hypergraph upper bound}. (Note that if one is not concerned about the dependence on $s$, then one can simply take $t=1$.)

\begin{theorem}[{\cite[Section~4]{ARS:99}}]\label{thm:ARS}
Given $r,t\in \bN$ and a prime $q\equiv 1\mod{t}$, there exists a $(q^{r-1}-1)$-regular $K_{r,s+1}$-free graph of order $q^{r-1}(q-1)/t$, where $s=(r-1)!t^{r-1}$.
\end{theorem}

In order to apply the above theorems, we also need the following fact about primes in arithmetic progressions, which follows from the Siegel--Walfisz theorem.

\begin{theorem}\label{thm:primes}
For every $\eps>0$ and $A\in \bN$, the following holds for sufficiently large~$x$. For any coprime integers $q,a$ with $q\le \log^A x$, there exists a prime $p\equiv a\mod{q}$ with $x\le p\le (1+\eps)x$.
\end{theorem}
\COMMENT{Siegel--Walfisz theorem asserts that the number of primes $\le x$ which are $\equiv a\mod q$ is $\pi(x;q,a)= Li(x)/\phi(q) \pm \cO_A(x\log^{-A} x)$. Since $Li(x)\sim x/\log x$ and $\phi(q)\le q$, we can deduce that $\pi((1+\eps)x;q,a) > \pi(x;q,a)$ when $q\le \log^{A-2} x$ and $x$ is sufficiently large.}

We are now ready to prove Conjecture~\ref{conj:yuster} and Theorem~\ref{thm:hypergraph upper bound}.

\lateproof{Conjecture~\ref{conj:yuster}}
Fix $\eps>0$ and suppose $k$ is sufficiently large. We aim to show that $f(k,2,s)\le (1+\eps)k^2 /s$ for every $1\le s\le k$.
The case when $s\gg \log k $ is already covered by Theorem~\ref{thm:hypergraph upper bound large s} proved above.
Thus, we may assume that $s\le \log^2 k$, say.
Let $p$ be the smallest prime at least $ k+2$ such that $p\equiv 1\mod{s}$. By Theorem~\ref{thm:primes}, we can assume $p\le (1+\eps/3)k\le \sqrt{1+\eps}k$. 
Let $G$ be the bipartite graph in Theorem~\ref{thm:mors}. Delete an arbitrary vertex from~$B$. Then every vertex in $A$ still has degree at least~$k$, and any two vertices in $A$ have at most $s$ common neighbours. Thus, $w(|A|,|A|-1;2,s)\ge k$, and Lemma~\ref{lem:Turan reduction} implies that $f(k,2,s)<|A|\le p^2/s \le (1+\eps)k^2 /s$, as desired.
\endproof

\lateproof{Theorem~\ref{thm:hypergraph upper bound}}
The lower bound follows from~\eqref{general LLL lower bound}. For the upper bound, we proceed similarly as in the proof of Conjecture~\ref{conj:yuster}. Since the case when $s\gg \log k $ is already covered by Theorem~\ref{thm:hypergraph upper bound large s},  we may assume that $(r-1)!\le s\le \log^2 k$, say, and that $k$ is sufficiently large.
Let $t$ be the maximal integer such that $(r-1)!t^{r-1}\le s$. Note that $t\ge 1$ and $t=\Theta_r(\sqrt[r-1]{s})$. Let $q$ be the smallest prime $\equiv 1\mod {t}$ with $q^{r-1}\ge 2rk$. Theorem~\ref{thm:primes} ensures that $q=\cO_r(\sqrt[r-1]{k})$. 
By Theorem~\ref{thm:ARS}, there exists a $(q^{r-1}-1)$-regular $K_{r,s+1}$-free graph $H$ of order $q^{r-1}(q-1)/t$. Now, partition $V(H)$ into parts $A$ and $B$ such that every vertex has at least $k$ neighbours in $B$ and $|A|>(r-1)|B|$. That this is possible can be seen via a simple probabilistic argument. For each vertex independently, put it into $B$ with probability $1/(r+1)$ and into $A$ otherwise. Then each vertex expects at least $(2rk-1)/(r+1)\ge 1.3k$ neighbours in~$B$ and $\expn{|A|}=r \cdot \expn{|B|}$. Using Chernoff’s bound, it is easy to see that the desired properties hold with high probability. 

Clearly, any $r$ vertices in $A$ have at most $s$ common neighbours. Hence, invoking Lemma~\ref{lem:Turan reduction}, we have $f(k,r,s)<|A|\le q^r/t =\cO_r(\sqrt[r-1]{k^r/s})$.
\endproof

\section{Concluding remarks}

\begin{itemize}
\item We studied sufficient conditions for the existence of independent transversals in multipartite graphs. Concretely, 
we have shown (Theorem~\ref{thm:main general}) that any multipartite graph with local degree $o(D)$, where each part has size at least $(1+o(1))D$ and average degree at most~$D$, contains an independent transversal. It would be very interesting to develop a hypergraph version of this result.

\item We completely settled the Erd\H{o}s--Gy\'{a}rf\'{a}s--{\L}uczak problem and its generalization due to Yuster in the graph case (Theorem~\ref{thm:yuster conj}). For hypergraphs ($r>2$), it would be desirable to prove that $f(k,r,1)=\cO_r(\sqrt[r-1]{k^r})$. We discovered a connection between upper bounds for $f(k,r,s)$ and lower bounds for $\ex(n;K_{r,s})$, thus indicating that the problem could be hard. On the other hand, our reduction does not work for $s<r-1$, so new ideas are needed for these cases, which could also shed new light on the general problem.


\item It could be illuminating to investigate random $(n,k,r,s)$-graphs further. The first moment method yields that if $n\gg \sqrt[r-1]{(k^r\log k)/s}$, then with high probability, a uniformly random $(n,k,r,s)$-graph has no independent transversal. What is the minimal $n$ such that a uniformly random $(n,k,r,s)$-graph has no independent transversal with probability at least~$1/2$? This seems open even for $r=2$ and $s=1$.
\end{itemize}

\section*{Acknowledgements}
We thank three anonymous reviewers for a careful reading of this paper and many helpful remarks.

\vspace{0.25cm}
\noindent
{\bf Note added in proof.} After this paper was released as preprint arXiv:2003.01683, Theorem~\ref{thm:main general} was proved independently by Kang and Kelly~\cite{KK:ta}, see arXiv:2003.05233. Whereas our main motivation came from the Erd\H{o}s--Gy\'{a}rf\'{a}s--{\L}uczak problem, they took the perspective of vertex colouring, and their proof is slightly different.

\providecommand{\bysame}{\leavevmode\hbox to3em{\hrulefill}\thinspace}
\providecommand{\MR}{\relax\ifhmode\unskip\space\fi MR }
\providecommand{\MRhref}[2]{%
  \href{http://www.ams.org/mathscinet-getitem?mr=#1}{#2}
}
\providecommand{\href}[2]{#2}

%

\begin{thebibliography}{10}

\bibitem{ABZ:07}
R.~Aharoni, E.~Berger, and R.~Ziv, \emph{Independent systems of representatives
  in weighted graphs}, Combinatorica~\textbf{27} (2007), 253--267.

\bibitem{alon:88}
N.~Alon, \emph{The linear arboricity of graphs}, Israel J. Math.~\textbf{62}
  (1988), 311--325.

\bibitem{alon:92}
\bysame, \emph{The strong chromatic number of a graph}, Random Structures
  Algorithms~\textbf{3} (1992), 1--7.

\bibitem{ARS:99}
N.~Alon, L.~R\'{o}nyai, and T.~Szab\'{o}, \emph{Norm-graphs: variations and
  applications}, J. Combin. Theory Ser. B~\textbf{76} (1999), 280--290.

\bibitem{AS:08}
N.~Alon and J.~H.~Spencer, \emph{The probabilistic method}, 3rd ed.,
  Wiley-Intersci. Ser. Discrete Math. Optim., John Wiley \& Sons, 2008.

\bibitem{BH:02}
T.~Bohman and R.~Holzman, \emph{On a list coloring conjecture of {R}eed}, J.
  Graph Theory~\textbf{41} (2002), 106--109.

\bibitem{BES:75}
B.~Bollob\'{a}s, P.~Erd\H{o}s, and E.~Szemer\'{e}di, \emph{On complete
  subgraphs of {$r$}-chromatic graphs}, Discrete Math.~\textbf{13} (1975),
  97--107.

\bibitem{brown:66}
W.~G.~Brown, \emph{On graphs that do not contain a {T}homsen graph}, Canad.
  Math. Bull.~\textbf{9} (1966), 281--285.
  
 \bibitem{erdos:72}
  P.~Erd\H{o}s, \emph{Problem~2}, in: Combinatorics, D.~J.~A.~Welsh and D.~R.~Woodall, eds. (The Institute of Mathematics and its Applications, 1972), 353--354.

\bibitem{EGL:94}
P.~Erd\H{o}s, A.~Gy\'{a}rf\'{a}s, and T.~{\L}uczak, \emph{Independent
  transversals in sparse partite hypergraphs}, Combin. Probab.
  Comput.~\textbf{3} (1994), 293--296.

\bibitem{ERS:66}
P.~Erd\H{o}s, A.~R\'{e}nyi, and V.~T.~S\'{o}s, \emph{On a problem of graph
  theory}, Studia Sci. Math. Hungar.~\textbf{1} (1966), 215--235.

\bibitem{furedi:96}
Zolt\'{a}n~F\"{u}redi, \emph{New asymptotics for bipartite {T}ur\'{a}n
  numbers}, J. Combin. Theory Ser. A~\textbf{75} (1996), 141--144.

\bibitem{haxell:16}
P.~Haxell, \emph{Independent transversals and hypergraph matchings---an
  elementary approach}, Recent trends in combinatorics, IMA Vol. Math. Appl.
  159, Springer, [Cham], 2016, pp.~215--233.

\bibitem{haxell:01}
P.~E.~Haxell, \emph{A note on vertex list colouring}, Combin. Probab.
  Comput.~\textbf{10} (2001), 345--347.

\bibitem{haxell:04}
\bysame, \emph{On the strong chromatic number}, Combin. Probab.
  Comput.~\textbf{13} (2004), 857--865.

\bibitem{JLR:00}
S.~Janson, T.~{\L}uczak, and A.~Ruci\'{n}ski, \emph{Random graphs},
  Wiley-Intersci. Ser. Discrete Math. Optim., Wiley-Interscience, 2000.

\bibitem{jin:92}
G.~P.~Jin, \emph{Complete subgraphs of {$r$}-partite graphs}, Combin. Probab.
  Comput.~\textbf{1} (1992), 241--250.
  
\bibitem{KK:ta}
 R.~J.~Kang and T.~Kelly, \emph{Colourings, transversals and local sparsity}, Random Structures
 Algorithms (to appear).

\bibitem{KRS:96}
J.~Koll\'{a}r, L.~R\'{o}nyai, and T.~Szab\'{o}, \emph{Norm-graphs and bipartite
  {T}ur\'{a}n numbers}, Combinatorica~\textbf{16} (1996), 399--406.

\bibitem{KST:54}
T.~K\H{o}v\'ari, V.~T.~S\'{o}s, and P.~Tur\'{a}n, \emph{On a problem of {K}.
  {Z}arankiewicz}, Colloq. Math.~\textbf{3} (1954), 50--57.

\bibitem{LohSudakov:07}
P.-S.~Loh and B.~Sudakov, \emph{Independent transversals in locally sparse
  graphs}, J. Combin. Theory Ser. B~\textbf{97} (2007), 904--918.

\bibitem{meshulam:03}
R.~Meshulam, \emph{Domination numbers and homology}, J. Combin. Theory Ser.
  A~\textbf{102} (2003), 321--330.

\bibitem{MR:02}
M.~Molloy and B.~Reed, \emph{Graph colouring and the probabilistic method},
  Algorithms and Combinatorics~23, Springer-Verlag, Berlin, 2002.

\bibitem{mors:81}
M.~M\"{o}rs, \emph{A new result on the problem of {Z}arankiewicz}, J. Combin.
  Theory Ser. A~\textbf{31} (1981), 126--130.

\bibitem{reed:99}
B.~Reed, \emph{The list colouring constants}, J. Graph Theory~\textbf{31}
  (1999), 149--153.

\bibitem{RS:02}
B.~Reed and B.~Sudakov, \emph{Asymptotically the list colouring constants are
  1}, J. Combin. Theory Ser. B~\textbf{86} (2002), 27--37.

\bibitem{ST:06}
T.~Szab\'{o} and G.~Tardos, \emph{Extremal problems for transversals in graphs
  with bounded degree}, Combinatorica~\textbf{26} (2006), 333--351.

\bibitem{yuster:97}
R.~Yuster, \emph{Independent transversals and independent coverings in sparse
  partite graphs}, Combin. Probab. Comput.~\textbf{6} (1997), 115--125.

\bibitem{yuster:97b}
\bysame, \emph{Independent transversals in {$r$}-partite graphs}, Discrete
  Math.~\textbf{176} (1997), 255--261.

\end{thebibliography}
%

\end{document}